\documentclass[12pt,a4paper]{article}
\usepackage{amsmath,amssymb,latexsym,amsthm}
\usepackage[english]{babel}
\usepackage[latin2]{inputenc}

\begin{document}

\newtheorem{Dfn}{Definition}
\newtheorem{Theo}{Theorem}
\newtheorem{Lemma}[Theo]{Lemma}
\newtheorem{Prop}[Theo]{Proposition}
\newtheorem{Coro}[Theo]{Corollary}
\newcommand{\Pro}{\noindent{\em Proof. }}
\newcommand{\Rem}{\noindent{\em Remark. }}

\title{Affine extensions of loops}
\author{\'Agota Figula and Karl Strambach}
\date{}
\maketitle
\footnotetext[1]{2000 {\em Mathematics Subject Classification:} 20N05, 51F25, 12D15, 51N30, 51M05.}
\footnotetext[2]{{\em Key words and phrases:}  loop, Bol loops, unitary and orthogonal geometries with positive index.}

\section{Introduction}

Most of the  known examples of loops $L$ with  strong relations to geometry have 
classical groups  as the groups generated by their left translations (\cite{Im}, 
\cite{Konrad}, \cite{Kolb},\cite{Gabrieli}, \cite{Kiechle}, Chapter 9, \cite{loops}, 
Chapters 22 and 25, \cite{figula}, \cite{figula2}). 
These groups $G$ may be seen as subgroups of the stabilizer of $0$ in the group of 
affinities of suitable affine spaces ${\cal A}_n$, and as the elements of the loops 
$L$ one can  often take certain projective subspaces of the hyperplane at infinity of 
${\cal A}_n$. The semidirect products $T \rtimes G$, where $T$ is the translation group 
of the affine space ${\cal A}_n$, have in many cases a geometric interpretation as motion 
groups of affine metric geometries. In the papers $\cite{figula}, \cite{figula2}$ three 
dimensional connected differentiable loops are constructed which have the connected 
component of the  motion group of the $3$-dimensional hyperbolic or pseudo-euclidean 
geometry as the group topologically generated by the left translations and which are Bol, 
Bruck or left A-loops. The set of the left translations of these loops induces on the 
plane at infinity the set of left translations of a loop isotopic to the hyperbolic plane 
loop (cf. \cite{loops}, Chapter 22, p. 280, \cite{Kolb}, p. 189). This and the fact that, 
up to our knowledge, there are only few known examples of sharply transitive sections in 
affine metric motion groups,  motivated us to seek a simple geometric procedure for an 
extension of a loop realized as the image $\Sigma ^{\ast }$ of a sharply transitive section 
in a subgroup $G^{\ast }$ of the projective linear group $PGL(n-1, \mathbb K)$ to a loop realized 
as the image of a sharply transitive section in a group $\Delta =T' \rtimes C$  of affinities of 
the $n$-dimensional space ${\cal A}_n=\mathbb K^n$ over a commutative field $\mathbb K$. Moreover, we 
desire that $T'$ is a large  subgroup of affine translations  and that $\alpha 
(C)=G^{\ast }$ holds for the canonical homomorphism $\alpha :GL(n,\mathbb K) \to PGL(n, \mathbb K)$. We 
show that this goal can be achieved if  in the $(n-1)$-dimensional 
projective hyperplane $E$ of infinity of ${\cal A}_n$  for $G^{\ast }$ there exists   an orbit ${\cal O}$ 
of $m$-dimensional subspaces such that $\Sigma ^{\ast }$ acts sharply transitively on 
${\cal O}$, if there is a subspace of dimension $(n-1-m)$ having empty intersection with any 
element of ${\cal O}$ and if the restriction of $\alpha ^{-1}$ to $\Sigma ^{\ast }$ defines 
a bijection from $\alpha ^{-1}(\Sigma ^{\ast })$ onto $\Sigma ^{\ast }$. 

In the third section we demonstrate that our construction successfully can be applied to 
sharply transitive sections in unitary and orthogonal groups $SU_{p_2}(n,F)$ of positive index $p_2$ over 
ordered pythagorean $n$-real fields $F$. In this way we obtain many non-isotopic topological loops.  The groups 
 generated by the left translations of these loops are  
 semidirect products $T \rtimes C$, where $T$ is the full translation group 
of 
${\cal A}_n$ and where $\alpha (C)$ is a non-solvable normal  subgroup of 
$\alpha (SU_{p_2}(n,F))$.

In the last section we take for the groups $G$ unitary or orthogonal Lie groups of 
any positive index in order to obtain differentiable loops $L$ such that the group
topologically generated by the left translations of $L$ is  a pseudo-unitary motion
group or the connected component of a
pseudo-euclidean motion group.

\section{Some basic notations of loop theory} 

A set $L$ with a binary operation $(x,y) \mapsto x \cdot y$ is called a loop 
if there exists an element $e \in L$ such that $x=e \cdot x=x \cdot e$ holds 
for all $x \in L$ and the equations $a \cdot y=b$ and $x \cdot a=b$ have 
precisely one solution which we denote by $y=a \backslash b$ and $x=b/a$. 
The left 
translation $\lambda _a: y \mapsto a \cdot y :L \to L$ is a bijection of $L$ 
for any $a \in L$. 
Two loops $(L_1, \cdot )$ and $(L_2, \ast )$ are isotopic if there are three bijections 
$\alpha , \beta , \gamma :L_1 \to L_2$ such that $\alpha (x) \ast \beta (y)=\gamma (x \cdot y)$ 
holds for any $x,y \in L_1$. 
A loop $(L, \cdot )$ is called topological if $L$ is a topological space and the mappings 
$(x,y) \mapsto x \cdot y$, $(x,y) \mapsto x \backslash y$, 
$(x,y) \mapsto y / x: L^2 \to L$ are continuous.  A loop $(L, \cdot )$ is called differentiable 
if $L$ is a $C^{\infty }$-differentiable manifold and the mappings 
$(x,y) \mapsto x \cdot y$, $(x,y) \mapsto x \backslash y$, 
$(x,y) \mapsto y / x: L^2 \to L$ are differentiable.   
\newline
A loop $L$ is a Bol loop if the identity $x(y \cdot xz)=(x \cdot yx)z$ holds. 
A Bruck loop is a Bol loop $(L, \cdot )$  satisfying the automorphic inverse property, i.e. 
the identity $(x \cdot y)^{-1}=x^{-1} \cdot y^{-1}$ for all $x,y \in L$. 
A loop $L$ is a left A-loop if each $\lambda_{x,y}= \lambda _{x y}^{-1} 
\lambda _x \lambda _y :L \to L$ is an automorphism of $L$. 
\newline
Let $G$ be the group generated by the left 
translations of $L$ and let $H$ be the stabilizer of $e \in L$ in the group 
$G$. 
The left translations of $L$ form a subset of $G$ acting on the cosets 
$\{x H; x \in G\}$ such that for any given cosets $aH$ and $bH$ there exists 
precisely one left translation $ \lambda _z$ with $ \lambda _z a H=b H$. 
\newline
Conversely let $G$ be a group, let H be a subgroup of $G$ and let $\sigma : G/H \to G$ be a 
section with $\sigma (H)=1 \in G$ such that the subset $\sigma (G/H)$ generates $G$ and 
acts sharply transitively on the space $G/H$ of the left 
cosets $\{x H, x \in G\}$ (cf. \cite{loops}, p. 18).  We call such a section sharply 
transitive. 
Then the multiplication  defined by 
$x H \ast y H=\sigma (x H) y H$ on the factor space $G/H$ or by $x \ast y=\sigma(xyH)$ 
on $\sigma (G/H)$ yields a loop $L(\sigma )$. If $N$ is the largest normal subgroup of 
$G$ contained in $H$ then the factor group $G/N$ is isomorphic to the group generated by 
the left translations of $L(\sigma )$. 
\newline
Two loops $L_1$ and $L_2$ having the same group $G$ as the group generated by the left 
translations 
and the same stabilizer $H$ of $e \in L_1, L_2$ are isomorphic if there is 
an automorphism of $G$ leaving $H$ invariant and mapping 
$\sigma_1(G/H)$ onto  $\sigma_2(G/H)$. 
The automorphisms of a loop $L$ corresponding to a sharply transitive section 
$\sigma :G/H \to G$ are given by the automorphisms of $G$ leaving $H$ and $\sigma (G/H)$ invariant.  
If two loops  are isotopic then the groups generated by their left translations are 
isomorphic (\cite{Pflugfelder}, Theorem III.2.7, p. 65).  
Loops $L$ and $L'$ 
 having the same group $G$ generated by their left translations  
 are isotopic if and only if there is a loop $L''$ 
isomorphic to $L'$ having $G$ again as the group generated by its left 
translations and there exists an inner automorphism $\tau $ of $G$ 
mapping  $\sigma ''(G/H)$ belonging to $L''$ onto the set 
$\sigma (G/H)$ corresponding to $L$ (cf. 
\cite{loops}, Theorem 1.11.  pp. 21-22).

\section{Affine extensions}

Let $G$ be a subgroup of the general linear group $GL(n, \mathbb K)$ over a commutative 
field $\mathbb K$. Denote by $\alpha $ the canonical epimorphism from $GL(n, \mathbb K)$ 
onto $PGL(n, \mathbb K)$. The kernel $Z$ of $\alpha $ is the centre of $GL(n, \mathbb K)$. 
Let $\tilde{H }$ be a subgroup of $G$ with  $Z \cap G \le {\tilde H}$ such that  for the pair 
$G^{\ast }=\alpha (G)$ and $H^{\ast }=\alpha ({\tilde H} )$ there exists a sharply 
transitive section $\sigma ^{\ast }: G^{\ast }/ 
H^{\ast } \to G^{\ast }$ determining a loop $L^{\ast }$.
Moreover, we assume that 
  $\Sigma ^{\ast }:=\sigma ^{\ast } (G^{\ast}/H^{\ast})$ generates  $G^{\ast }$ and that  
 for the preimage $(\alpha |G)^{-1}(\Sigma ^{\ast })=\Sigma \subseteq  G$  
 one has ${\tilde H} \cap \Sigma= \{ 1 \}$. Then the mapping $\alpha $ induces a bijection from 
$\Sigma $ onto $\Sigma ^{\ast }$.

We denote by ${\cal A}_n$ the $n$-dimensional affine space $\mathbb K^n$ and by $E$ the projective hyperplane of dimension $(n-1)$ at infinity of ${\cal A}_n$. Let $U^{\ast }$ be an $m$-dimensional subspace of $E$ having  $ H^{\ast }$ as the stabilizer of $U^{\ast }$ in 
$G^{\ast }$. Let ${\cal X} $ be the set  
\[ {\cal X} = \{\gamma  U^{\ast }; \gamma \in \Sigma ^{\ast } \}.\]
The elements of ${\cal X} $ may be seen as the elements of $L^{\ast }$ such that  
$U^{\ast }$  is the identity of $L^{\ast }$ and the multiplication is given by 
$X^{\ast } \circ Y^{\ast } = \tau ^{\ast }_{U^{\ast },X^{\ast }}(Y^{\ast })$ for all 
$X^{\ast },Y^{\ast } \in {\cal X} $, where $\tau^{\ast }_{U^{\ast },X^{\ast }}$ is the 
unique element of the sharply transitive set $\Sigma ^{\ast }$ of the linear transformations 
of $E$ mapping $U^{\ast }$ onto $X^{\ast }$. 

Let $A=T \rtimes S $ be the semidirect product consisting of  affinities of  
${\cal A}_n=\mathbb K^n$, where $T$ is the translation group of ${\cal A}_n$ and $S$ is the stabilizer of  
$0 \in {\cal A}_n$  isomorphic to the group $GL(n, \mathbb K)$.
 We consider the group $G$ as a subgroup of $S$ in the group 
$\Theta =\mathbb K^n \rtimes G $ of affinities of ${\cal A}_n$. 
 The subgroup ${\tilde H }$  of $S$  fixes the point $0 \in {\cal A}_n$ and the subspace $U^{\ast }$ of the hyperplane $E$.  Let $U$ be 
the $(m+1)$-dimensional affine subspace containing $0$ and intersecting $E$ in $U^{\ast }$.  
If  $H$ is the stabilizer of $U$ in the group $\Theta $, then  one has 
${\tilde H} =H \cap \Theta _0$, where $\Theta _0$ is the stabilizer of the point $0$ in 
$\Theta $. 

Let $W$ be a subspace of  ${\cal A}_n$ such that $W$ contains $0$, has affine dimension $(n-m-1)$ and 
intersects any subspace of the set ${\cal Z}:=\{ \rho (U); \rho \in \Sigma \}$  only in $0$. Let $T_W$ be the group 
of 
affine translations $x \mapsto x+w: {\cal A}_n \to {\cal A}_n$ with $w \in W$. Then 
$W$ intersects any subspace $\delta (Y)$, where $\delta \in T_W$ and $Y \in {\cal Z}, $ 
in precisely one point. Moreover, the stabilizer of  $\delta (Y)$ in $T_W$ 
consists only of the identity.    

\begin{Theo} The subset $\Xi = T_W \Sigma = \{ \tau \rho; \tau \in T_W, \rho 
\in \Sigma \}$ 
of the group $\Theta = T \rtimes G$ acts sharply transitively on the set 
\[ {\cal U} = \{ \psi (U); \psi \in \Xi \}=\{ \psi (U); \psi \in \Theta \} . \]
The elements of ${\cal U} $ can be taken as the elements of a loop $L_{\Xi }$ which 
has $U$ as the identity and for which the multiplication is defined by 
\[ X \circ Y=\tau _{U,X} (Y) \quad \hbox{for  all} \quad X,Y \in {\cal U} , \]
where $\tau _{U,X}$ is the unique element of $\Xi  $  mapping $U$ onto $X$. 

The set $\Xi $ is the set of the left translations of $L_{\Xi }$ and generates a group  $\Delta $ which is a semidirect product  $\Delta =T' \rtimes C$,
where the normal subgroup $T'$ consists of translations of the affine space ${\cal A}_n$ and  $C$ is a subgroup of $G$ with $\alpha (C)=G^*$.  

There is a sharply transitive section $\sigma : \Delta / {\hat H} \to \Delta$ such that  $\sigma(\Delta / {\hat H})= \Xi $, the group  ${\hat H}$ is the stabilizer  of $U$ in  $\Delta $ and the subgroup $T' \cap {\hat H}$ consists of all translations $x \mapsto x+u: {\cal A}_n \to {\cal A}_n$ with $u \in U$. 
\end{Theo}

\Pro Let $D_1$ and $D_2$ be elements belonging to ${\cal U} $. We show that there is 
precisely one element $\beta \in \Xi$ with $\beta (D_1)=D_2$. Let $D^{\ast }_1=D_1 \cap E$ 
and  $D^{\ast }_2=D_2 \cap E$, where $E$ is the hyperplane at infinity of ${\cal A}_n$. 
Thus  there exists precisely one element $\rho ^{\ast } \in \Sigma ^{\ast }$ and hence there exists 
precisely one element $\rho \in \Sigma $ with $\alpha (\rho )= \rho ^{\ast }$ such that 
$\rho ^{\ast } (D^{\ast }_1)=D^{\ast }_2$. The subspaces $\rho (D_1)$ and $D_2$ 
intersect $E$ in $D^{\ast }_2$. In the group $T_W$ there exists precisely one translation 
$\tau $ mapping the point $\rho(D_1) \cap W$ onto the point $D_2 \cap W$. Hence the 
element $\beta =\tau \rho$ is the only element in $\Xi $ 
mapping $D_1$ onto $D_2$ and the set $\Xi $ is a sharply transitive set on ${\cal U} $. 
It follows that the subspaces in ${\cal U}$ can be taken as the elements of a loop 
$L_{\Xi }$ having $U$ as the identity, such that the multiplication is defined as in 
the assertion of the theorem. 

The group $\Delta $ generated by the left translations of $L_{\Xi }$ is a subgroup of $\Theta = T \rtimes G$. Let ${\hat H}$ be the stabilizer of $U$ in  $\Delta $. Since $\Xi $ is the image of a sharply transitive section 
$\sigma :\Delta /{\hat H} \to \Delta$ we have $\Delta (U) = \Xi {\hat H}(U)= \Xi (U)$.   
Let $T_U$ be the group of affine translations 
$ x \mapsto x+u: {\cal A}_n \to {\cal A}_n$ with $u \in U $.
Since $W \oplus U=\mathbb K^n$ we have that $T=T_W \times T_U$.   
Thus  one has $\Delta T(U)=\Delta T_W T_U (U)=\Delta T_W (U)=\Delta (U)$ since $T_W \le  \Delta$. For the group $\Lambda $ of dilatations $x \mapsto a x:{\cal A}_n \to {\cal A}_n$ with $a \in  \mathbb K \backslash \{ 0 \}$ we have 
that $T \Lambda $ 
is a normal subgroup of $\Theta \Lambda $ and $\Lambda (U)=U$. Moreover 
$\Theta (U)= \Delta T \Lambda (U)$ since the kernel of the restriction of  
$\alpha :GL(n, \mathbb K) \to PGL(n, \mathbb K)$ to $G$ consists only of dilatations.  

The group $\Delta $ contains a normal subgroup $N$ fixing the hyperplane $E$ at  infinity pointwise. Since $\Sigma ^{\ast }$ generates $G^{\ast }$ we see that  $\Delta / N$ is isomorphic to $G^{\ast }$. 

Let $T'=T \cap \Delta $. Then  $\Delta $ is the semidirect product of $\Delta = T' \rtimes C$, where $C$ is the stabilizer of $0$ in $\Delta $ and 
$CN/N$ is isomorphic to $G^{\ast }$. 
\qed 

\section{Applications} 
Let $R$ be an ordered pythagorean field and let $K=R(i)$ be the algebraic extension of $R$ 
such that $i^2=-1$. 
 Let $F \in \{ R, K \}$ and let $V=F^n$ be an $n$-dimensional $F$-vector space for a fixed 
 $n \ge 3$. The automorphism  $a \mapsto {\bar a}: F \to F$ is the identity if $F=R$ or the involutory 
 automorphism fixing $R$ elementwise and mapping $i$ onto $-i$ if $F=K$. Denote by 
 ${\cal M}_n(F)$ the set of the $(n \times n)$-matrices over $F$. 
If $A=(a_{i,j})$ is a matrix in ${\cal M}_n(F)$ then $\bar{A}^t=({\bar a}_{j,i})$. Let 
${\cal H}(n,F)$ be the set of positive definite hermitian $(n \times n)$-matrices, i.e. 
the set \[ {\cal H}(n,F)=\{ A \in {\cal M}_n(F); A= \bar{A}^t \ \hbox{with } \ {\bar v}^t 
A v > 0 \ \hbox{for all} \ v \in V \backslash \{ 0 \} \} .\] 
We assume that the field $R$ is $n$-real which means that the characteristic polynomial of 
every matrix in ${\cal H}(n,F)$ splits over $K$ into linear factors. Thus  this polynomial  splits  into linear factors
already over $R$ (cf. \cite{Kiechle}, p. 14). The class of 
$n$-real fields contains the class of totally real fields (cf. \cite{Kiechle}, p. 13), 
which is larger than the class of real closed fields and the class of hereditary euclidean fields. A 
hereditary euclidean field $k$ is an ordered field such that every formally real algebraic 
extension of $k$ has odd degree over $k$ (cf. \cite{Prestel}, Satz 1.2 (3), p. 197). 

The group 
\[ U(n, F)=\{ B \in GL(n,F); B {\bar B}^t=I_n \} ,\] where $I_n$ is the identity in  $GL(n,F)$, is called the orthogonal group for 
$F=R$ and the unitary group for $F=K$. Let 
\[ J_{(p_1,p_2)}=\hbox{diag} (1, \ldots ,1,-1, \ldots , -1) \] be the diagonal 
$(n \times n)$-matrix such that the first $p_1$ entries are $1$ and the remaining $p_2$ 
entries are $-1$. We have $p_1+p_2=n$. The matrix $J_{(p_1,p_2)}$ defines 
a hermitian form on $F^n$ for $F=K$  
and an orthogonal form for $F=R$   by 
\[ {\bar v}^t J v= \sum ^{p_1}_{i=1} {\bar v}_i v_i - \sum ^n_{j=p_1+1} {\bar v}_j v_j . \]
Let $p_2 >0$. The unitary (orthogonal) group of index $p_2$ is the set 
\[ U_{p_2}(n,F)= \{ A \in GL_n( F); {\bar A}^t J_{(p_1,p_2)} A= J_{(p_1,p_2)} \} . \]
Since  the 
group $U_{p_2}(n,F)$ is isomorphic to the group $U_{(n-p_2)}(n,F)$ (cf. \cite{Porteous}, 
Proposition 9.11, p. 153)  we may assume that $p_1 \ge p_2$. 
Let \[ \Omega _{(p_1,p_2)}(F)= U_{p_2}(n,F) \cap U(n, F) \  \ \hbox{and}
\ \  \Sigma _{(p_1,p_2)}(F)= U_{p_2}(n,F) \cap {\cal H}(n, F). \]
The group $\Omega _{(p_1,p_2)}(F)$ is the direct product $\Omega _{(p_1,p_2)}(F)=U(p_1, F) 
\times U(p_2,  F)$, where 
$U(p_1, F)$ may be identified with the group $\left ( \begin{array}{cc}
U(p_1,  F) & 0 \\
0 &  I_{p_2} \end{array} \right )$ and $U(p_2, F)$ may be identified with the group 
$\left ( \begin{array}{cc}
 I_{p_1}  & 0 \\
0 &  U(p_2, F) \end{array} \right )$; here  $I_{p_i}$ is the identity in $GL(p_i, F)$ 
(cf. \cite{Kiechle}, Theorem 9.13, p. 123). 

According to \cite{Kiechle} (Theorem 9.11, p. 121) the set 
$\Sigma _{(p_1,p_2)}(F)$  is the image of a sharply transitive section 
$\sigma ' : U_{p_2}(n,F)/ \Omega _{(p_1,p_2)}(F) \to U_{p_2}(n,F)$ such that the corresponding 
loop $L_{(p_1,p_2)}$ is a Bruck loop.  

The group  $G_{(p_1,p_2)}$ generated by the set  $\Sigma _{(p_1,p_2)}(F)$ of the left translations of  $L_{(p_1,p_2)}$ is contained in the group  $SU_{p_2}(n,F):=\{ A \in  U_{p_2}(n,F);\ \hbox{det}\ A=1\}$ (cf. \cite{Kiechle}, 9.14, p. 124). Thus the loop $L_{(p_1,p_2)}$ corresponds also to the section 
\newline
\centerline{$\sigma: SU_{p_2}(n,F)/ \Phi \to SU_{p_2}(n,F)$,} 
where $\Phi:=(U(p_1, F) \times U(p_2,  F)) \cap  SU_{p_2}(n,F)$. 

The kernel of the restriction of $\alpha : GL(n, F) \to PGL(n,F)$ to the group $SU_{p_2}(n,F)$ consists of the matrices  
$D_a=\hbox{diag} (a, \ldots , a), \ a \in F \backslash \{ 0 \}$ and $a^n=1$. Moreover one has  $a {\bar a}=1$ since any matrix $D_a$ satisfies ${\bar D_a}^t J_{(p_1,p_2)} D_a= J_{(p_1,p_2)}$. Thus  any matrix $D_a$ is contained in  $\Phi $ and $\alpha $ induces a bijection from  $\Sigma _{(p_1,p_2)}(F)$ onto 
$\alpha ( \Sigma _{(p_1,p_2)}(F))$. The set 
$\alpha ( \Sigma _{(p_1,p_2)}(F))$ is the image of a sharply transitive 
section 
\newline
\centerline{$\sigma ^*:\alpha (SU_{p_2}(n,F))/\alpha (\Phi ) \to \alpha (SU_{p_2}(n,F)) $}
 which corresponds to a Bruck loop $L^*_{(p_1,p_2)}$.

The elements of $\Sigma _{(p_1,p_2)}(F)$ are matrices $A \in  SU_{p_2}(n,F)$ satisfying the relations $A= \bar{A}^t$ and ${\bar v}^t 
A v > 0$ for all $v \in V \backslash \{ 0 \}$. With $A$ also $A^{-1}$ is contained in  $\Sigma _{(p_1,p_2)}(F)$ (\cite{Kiechle} 1.11, p. 16).
Because of $B^{-1}={\bar B}^t$ for all  $B \in \Phi$ 
and 
${\bar B}^t A B \in  \Sigma _{(p_1,p_2)}(F)$ (\cite{Kiechle} 1.11, p. 16)
one has 
\begin{equation}  B^{-1} A B \in \Sigma _{(p_1,p_2)}(F) \ \ \hbox{for \ all}\  B \in \Phi \  \hbox{and}\  A \in \Sigma _{(p_1,p_2)}(F).   \end{equation} 
   
Since $\sigma $ is a section every element $S$ of  $SU_{p_2}(n,F)$ can be written in a unique way as $S=S_1 C$ with $S_1 \in  \Sigma _{(p_1,p_2)}(F)$ and $C \in \Phi$. The set 
\[ \Sigma _{(p_1,p_2)}(F)^{G_{(p_1,p_2)}}=\{Y^{-1} X Y;\ X \in \Sigma _{(p_1,p_2)}(F), Y \in G_{(p_1,p_2)} \} \]
is invariant with respect to the conjugation by the elements 
$S \in SU_{p_2}(n,F):$
\[ S^{-1} Y^{-1} X Y S=C^{-1} S_1^{-1} Y^{-1} X Y S_1 C= \]
\[[(C^{-1} S_1^{-1} C)(C^{-1} Y^{-1} C)](C^{-1} X C)[(C^{-1} Y C)(C^{-1} S_1 C)]\in \Sigma _{(p_1,p_2)}(F)^{G_{(p_1,p_2)}}. \]
Hence the group $G_{(p_1,p_2)}$, which is generated also by 
$\Sigma _{(p_1,p_2)}(F)^{G_{(p_1,p_2)}}$, is a normal non central subgroup of $SU_{p_2}(n,F)$. Then 
according to  Th\'eor\'eme 5 in \cite{Dieudonne} p. 70  the group $G_{(p_1,p_2)}$ coincides with $SU_{p_2}(n,F)$ if $F=K$. If $F=R$ and $(n,p_2) \neq (4,2)$ then the group   $G_{(p_1,p_2)}$ contains the commutator subgroup 
$[SU_{p_2}(n,F)]'=:{\mathcal K}_{(n,p_2)}$ of $SU_{p_2}(n,F)$ (\cite{Dieudonne2}, p. 63 and pp. 58-59). 
If $F=R$ and $(n,p_2)= (4,2)$ then the commutator subgroup ${\mathcal K}_{(4,2)}$ is isomorphic to the direct product $PSL_2(R) \times PSL_2(R)$ (\cite{Dieudonne2}, p. 59). Since the hermitian matrices in the set $\Sigma _{(2,2)}(F)$ depend on $3$ free parameters (\cite{Kiechle}, 9.12, p. 122) the group $G_{(2,2)}$ contains ${\mathcal K}_{(4,2)}$. Therefore 
in any case the group $G_{(p_1,p_2)}$ is a normal subgroup  
of $SU_{p_2}(n,F)$ containing ${\mathcal K}_{(n,p_2)}$.

The group $G_{(p_1,p_2)}$ 
leaves the value 
${\bar v}^t J_{(p_1,p_2)} v$ invariant since 
\[{\bar v}^t ({\bar A}^t J_{(p_1,p_2)} A) v= {\bar v}^t J_{(p_1,p_2)} v\ \hbox{ for}\  
A \in SU_{p_2}(n,F). \] 
We  see the group $G_{(p_1,p_2)}$ as a subgroup of the stabilizer of the element $0$ 
in the group $A$ of affinities  of  ${\cal A }_n=F^n$, and the group $\alpha (G_{(p_1,p_2)}):=G^*_{(p_1,p_2)}$  
 as a subgroup of the group $PGL(n,F)$ which acts on the $(n-1)$-dimensional  projective hyperplane $E$ at infinity of ${\cal A }_n$. 

We embed the affine space ${\cal A }_n$ into the $n$-dimensional projective space $P_n(F)$ 
such that $(x_1,\cdots ,x_n) \mapsto  
F^{\ast } (1,x_1,\cdots ,x_n)$, $x_i \in F$ for all $1 \le i \le n$ and  
$F^{\ast }= F \backslash \{ 0 \}$. With respect to 
this embedding the hyperplane $E$ consists  of the points $\{ F^{\ast } (0,x_1,\cdots ,x_n), x_i \in F, \hbox{not 
all } x_i=0 \}$. The cone in ${\cal A }_n$ which is described by the equation 
\newline
\newline
$(\ast )$ \centerline{ $\sum \limits_{i=1}^{p_1} {\bar x}_i x_i - 
\sum \limits_{j=p_1+1}^{n} {\bar x}_j x_j=0$ }
\newline
\newline
intersects $E$ in a hyperquadric $C$; the points 
$\{ F^{\ast } (0,x_1,\cdots ,x_n)\}$ of $C$ satisfy the equation $(\ast )$. 
The hypersurface $C$ of $E$ divides the points of $E \backslash C$ into 
two regions $R_1$ and $R_2$. A point $F^{\ast } (0,x_1,\cdots ,x_n)$ belongs to $R_1$ if 
and only if $\sum \limits_{i=1}^{p_1} {\bar x}_i x_i > \sum \limits_{j=p_1+1}^{n} {\bar x}_j x_j $.
It belongs to $R_2$ if and only if 
$\sum \limits_{i=1}^{p_1} {\bar x}_i x_i < \sum \limits_{j=p_1+1}^{n} {\bar x}_j x_j $. 
The group $\alpha(SU_{p_2}(n,F))=SU_{p_2}(n,F)/\Lambda '$, where $\Lambda ' $ is the group of dilatations contained in $SU_{p_2}(n,F)$, 
leaves $R_1$, $R_2$ as well as $C$ invariant 
since for any $f \in F$ and $v \in V=F^n$ one has 
$({\bar f}{\bar v}^t) J_{(p_1,p_2)} (f v)= ({\bar f} f)({\bar v}^t J_{(p_1,p_2)} v)$ and ${\bar f} f
>0$. The group $\alpha (\Phi )= \Phi /( \Phi \cap \Lambda ' )$ 
 leaves the subspace 
\[ W_1^{\ast }=\{ (0,x_1,\ldots ,x_{p_1},0, \ldots ,0); x_i \in F \} \subseteq E\] 
as well as the subspace 
\[ W_2^{\ast }=\{ (0,\ldots ,0,x_{p_1+1}, \ldots ,x_n); x_i \in F \} \subseteq E  \]  
invariant. The intersection  $W_1^{\ast } \cap W_2^{\ast }$ is empty since $ 
W_i^{\ast } \subseteq R_i, i=1,2$. 

Let $W_i, i=1,2$, be the $p_i$-dimensional 
affine subspace of ${\cal A}_n$ containing $0$ such that $W_i \cap E=W_i^{\ast }$. 
Thus  $W_1 \cap W_2=\{ 0 \}$. Let $\tilde{W}_j$ be a $p_j$-dimensional affine subspace 
of ${\cal A}_n$ such that $p_j=n-p_i$ and  $\tilde{W}_j$ intersects $W_i$ only in 
the point $0$. Thus $\tilde{W}_j$ intersects any subspace of the set 
\[ {\cal Z}_i= \{ \rho (W_i), \rho \in G_{(p_1,p_2)} \}=\{ \lambda (W_i), \lambda \in 
\Sigma _{(p_1,p_2)}(F) \}, \] where $i \neq j$, only in $0$.  Affine 
subspaces $\tilde{W}_j$ with these properties exist, one can take for instance 
$\tilde{W}_j=\rho (W_j) \in  {\cal Z}_j$.  

Let $\Theta $ be the 
semidirect product $\Theta =T \rtimes G_{(p_1,p_2)}$, where $T$ is the translation 
group of ${\cal A}_n$.  By Theorem 1  the set 
$\Xi_{(p_i, {\tilde W}_j)}=\{ T_{{\tilde W}_j} \Sigma _{(p_1,p_2)}(F) \}$, $i \neq j$, 
acts sharply transitively on the set 
\[ {\cal U}_i= \{ \psi (W_i); \psi \in \Xi_{(p_i, {\tilde W}_j)} \}.  \]
Thus  a loop $L_{(p_i, {\tilde W}_j)} $ is realized on  ${\cal U}_i$. 

The group $SU_{p_2}(n,K)$ acts irreducibly on the vector space $V=K^n$ and the commutator subgroup ${\mathcal K}_{(n,p_2)}$ of  $SO_{p_2}(n,R)$ acts irreducibly on $V=R^n$ (cf. \cite{Artin}, Theorem 3.24, p. 136). Hence the group 
 $\Delta $ generated by the left translations 
$ \Xi_{(p_i, {\tilde W}_j)}$ of the loop $L_{(p_i, {\tilde W}_j)} $ contains  all translations of the affine space ${\cal A}_n$. It follows that  $\Delta $ is  the semidirect product 
$\Delta =T \rtimes C$ of the translation group $T$ by a  
 subgroup $C$ of the stabilizer  of $0 \in {\cal A}_n$ 
in the group $A$ of affinities. If $F=K$ then $C$ is isomorphic to $SU_{p_2}(n,K)$ and the stabilizer $\hat{H}$ of 
$W_i$ in $\Delta $ is the semidirect product $T_{W_i} \rtimes \Phi$  since any element 
$g \in  G_{(p_1,p_2)}=SU_{p_2}(n,K)$ has a unique representation as $g=g_1 g_2$ with 
$g_1 \in \Sigma _{(p_1,p_2)}(K)$ and $g_2 \in \Phi$. If $F=R$ then $C$ is a normal subgroup of $SO_{p_2}(n,R)$ 
containing  ${\mathcal K}_{(n,p_2)}$ and the stabilizer $\hat{H}$ of 
$W_i$ in $\Delta $ is the semidirect product $T_{W_i} \rtimes \Gamma $, where $\Gamma =\Phi \cap C$.

{\bf For $p_1 > p_2$ the loop $L_{(p_1, {\tilde W}_2)} $ is  never isotopic to a loop 
$L_{(p_2, {\tilde W}_1)} $.} This follows from the fact that the stabilizer $H_k$, 
$k=1,2$, of the identity of $L_{(p_k, {\tilde W}_l)} $ with $l \neq k$ in $\Delta $ 
contains the group $T_{W_k}$ as the largest normal subgroup consisting of affine translations. 
 Since $T_{W_1}$ is not isomorphic to $T_{W_2}$ one has that $H_1$ is not isomorphic to $H_2$. (cf. \cite{Pflugfelder}, Theorem III.2.7, p. 65)

Now we consider 
 the loops $L_{(p_i, W_j)}$ and $L_{(p_i, {\tilde W}_j)}$ for $W_j \neq  {\tilde W}_j$. According to $(1)$ 
 the subspaces $W_1$ and $W_2$ are invariant under the subgroup $\Phi $ 
 of the stabilizer of  $0 \in {\cal A}_n$  in the group $A$ of affinities. Hence 
 if $g \in \Phi $ then  one has $g \Sigma _{(p_1,p_2)}(F) g^{-1}=\Sigma 
 _{(p_1,p_2)}(F)$ and $g T_{W_k} g^{-1}=T _{W_k}$, $k=1,2$, for the group 
$T_{W_k}=\{ x \mapsto x+w_k; w_k \in W_k \}$.  
 This yields $g 
\Xi _{(p_i, W_j)} g^{-1} =\Xi _{(p_i, W_j)}$. For 
 $W_j \neq  {\tilde W}_j$ the group $\Phi $
does not normalize the translation group $T_{{\tilde W}_j}$. Therefore 
\[ g T_{{\tilde W}_j} \Sigma _{(p_1,p_2)}(F) g^{-1}= (g T_{{\tilde W}_j} g^{-1}) (g 
 \Sigma _{(p_1,p_2)}(F) g^{-1})= \] \[(g T_{{\tilde W}_j} g^{-1})  \Sigma _{(p_1,p_2)}(F) 
 \neq \Xi _{(p_1, {\tilde W}_j)} \] 
 for suitable elements $g \in \Phi $. This means that not all elements 
 of $ \Phi $ induce automorphisms of $L_{(p_i, {\tilde W}_j)}$. 
Therefore {\bf the loops $L_{(p_i, W_j)}$ and $L_{(p_i, {\tilde W}_j)}$ are not 
isomorphic if $W_j \neq  {\tilde W}_j$.}  

\begin{Prop} Any loop $L_{(p_i, {\tilde W}_j)}$ is a topological loop with respect 
to the topology induced on the set ${\cal U} $ by the topology on the set of the 
$p_i$-dimensional subspaces of ${\cal A}_n$ which is derived from the topology of 
the topological field $F$.  
\end{Prop}
\Pro
Since $R$ is an ordered field, $R$ as well as $K=R(i)$ are topological fields with respect 
to the topology given by the ordering of $R$. Then the ring 
${\cal M}_n(F)$ of $(n \times n)$-matrices over $F$ is a topological ring such that the open 
$\varepsilon $-neighbourhoods of $0 \in {\cal M}_n(F)$ consist of matrices $(c_{i,j})$ with 
$|c_{ij} | < \varepsilon $. The group $GL(n, F) \le {\cal M}_n(F)$ is a topological group. 
Since the set  $Z= \{ \hbox{diag} (a, \ldots , a), a \in F \backslash \{ 0 \} \} $ is a 
closed subgroup of $GL(n, F)$ the group $PGL(n,F)= GL(n,F)/ Z$ is a topological group. 
The subgroups $SU_{p_2}(n,F)$ 
and $\Phi=(U(p_1, F) \times U(p_2,  F)) \cap  SU_{p_2}(n,F)$ are closed subgroups of $GL(n,F)$. Moreover $SU_{p_2}(n,F)Z /Z$ as 
well as  $\Phi Z/Z$ are closed subgroups of $PGL(n,F)$. 

The affine space 
${\cal A}_n=F^n$ and the $(n-1)$-dimensional 
projective  hyperplane $E$  carry topologies derived from the topology of the field $F$ 
(cf. \cite{Lenz}, Chapter XI). The semidirect product $A=T \rtimes GL(n,F)$ 
is a topological group consisting of continuous affinities; it  induces on the hyperplane $E$ a continuous group of projective collineations. Any subset of $A$ is a topological space with respect to the topology induced from $A$ and any subgroup of $A$ becomes  a topological group in this manner.

Let $Q_1$ be a fixed 
$p_i$-dimensional subspace of  ${\cal A}_n$  and let ${\cal Q}$ 
be the set of the affine $(n-p_i)$-dimensional affine subspaces with $| Q_1 \cap Q|=1$ for 
$Q \in {\cal Q}$. The set ${\cal Q}$ also carries  a topology  determined by the 
topology of $F$. The set ${\cal Q}^{\ast }$ of intersections  
$Q^{\ast }$ of the affine subspaces $Q$ of ${\cal Q}$ with $E$ inherits the topology of 
the Gra\ss mannian manifold of the $(n-p_i-1)$-dimensional subspaces of the hyperplane $E$. The geometric operation 
$ (Q, Q_1) \mapsto Q \cap Q_1 : {\cal Q} \to Q_1 $ is continuous. 

On the topological space $\Sigma _{(p_1,p_2)}(F)$ a topological Bruck loop $L_{(p_1,p_2)}$ is realized by 
 the multiplication  
\begin{equation} A \circ B= \sqrt{A B^2 A} \ \hbox{for all } \ A,B \in 
\Sigma _{(p_1,p_2)}(F),  \end{equation}
where $X \mapsto \sqrt{X}$  
is the inverse map of the bijection $X \mapsto X^2: \Sigma _{(p_1,p_2)}(F) \to 
\Sigma _{(p_1,p_2)}(F)$  
(cf. \cite{Kiechle} (1.14), p. 17 and (9.1) Theorem (4), p. 108, \cite{loops}, p. 121). 
\newline
We denote by $[\rho (W_i)]^{\ast }$ with $\rho \in \Sigma _{(p_1,p_2)}(F)$
 the 
intersection of the  subspace $\rho (W_i)$ with the hyperplane $E$ and by ${\cal Z}_i^*$ the set $\{ [\rho (W_i)]^{\ast };\ \rho \in  \Sigma _{(p_1,p_2)}(F) \}$. For the elements of the loop $L_{(p_i, {\tilde W}_j)}$ one can take the elements  of the set 
\[ {\mathcal U}_{(p_i, {\tilde W}_j)}=\{ \psi (W_i);\ \psi \in \Xi_{(p_i, {\tilde W}_j)} \}=\{ \tau \rho (W_i);\ \tau \in T_{{\tilde W}_j}, \rho \in  \Sigma _{(p_1,p_2)}(F) \}.\]  
The subspace ${\tilde W}_j$ is homeomorphic to the group $T_{{\tilde W}_j}$, 
and the set ${\cal Z}_i$ is homeomorphic to $\Sigma _{(p_1,p_2)}(F)$. 
Any element $ \tau \rho (W_i) \in {\mathcal U}_{(p_i, {\tilde W}_j)}$ is uniquely determined by $[\rho (W_i)]^{\ast }$ and  $(\tau \rho (W_i)) \cap {\tilde W}_j$. The mapping 
\[ \omega : \tau \rho (W_i) \mapsto ((\tau \rho (W_i)) \cap {\tilde W}_j,[\rho (W_i)]^{\ast }) \]
from ${\mathcal U}_{(p_i, {\tilde W}_j)}$ onto the topological product $ {\tilde W}_j \times {\cal Z}_i^*$ is a bijection such that  
\[ \omega^{-1}:(w, Z^*) \mapsto w \vee Z^*, \]
where $w \vee Z^*$ is the $p_i$-dimensional affine subspace containing $w \in  {\tilde W}_j$ 
and intersecting  $E$ in $Z^* \in  {\cal Z}_i^*$. 
Since the geometric operations of 
joining and of intersecting of distinct subspaces are continuous maps,  $\omega $ is a homeomorphism.  

Let $(w_k, Z_k^*) \in {\tilde W}_j \times {\cal Z}_i^*$ with $k=1,2$ and let  
$ \tau_k \rho_k (W_i)$ be the subspaces of  ${\mathcal U}_{(p_i, {\tilde W}_j)}$ such that $\omega(\tau_k \rho_k (W_i))=(w_k, Z_k^*)$.  
The 
multiplication given by 
\begin{equation} (w_1, Z_1^*) \circ (w_2, Z_2^*)=(w_3, Z_3^*),   \end{equation} 
where $Z_3^*=[\rho_1 \rho_2(W_i)]^*$  and 
\[ w_3=\tau_1[\rho_1 (\tau_2 \rho_2(W_i)) \cap {\tilde W}_j]= \tau _1 \big[ \big(\rho _1[ \tau _2 \rho_2(W_i) \cap {\tilde W}_j] \vee [\rho_1 \rho_2(W_i)]^{\ast } \big) \cap {\tilde W}_j \big] \] yields a topological loop. This loop is homeomorphic to  $L_{(p_i, {\tilde W}_j)}$ since 
\newline
$[\rho_1 \tau_2 \rho_2(W_i)]^*=[\rho_1 \rho_2(W_i)]^*$ and 
$[\tau_1 \rho_1 \tau_2 \rho_2(W_i)]^*=[\rho_1 \rho_2(W_i)]^*$.      
\qed 

\section{Special cases:  $\mathbb R$ and  $\mathbb C$}

\begin{Prop} The loop $L_{(p_i, {\tilde W}_j)}$ is a differentiable loop diffeomorphic to
$\mathbb R^d$, where $d = \varepsilon (p_j+p_1 p_2)$, with  $\varepsilon =1$ if $F = \mathbb R$ and $\varepsilon =2$ 
if $F= \mathbb C$.  

If  $F= \mathbb C$ then the group $\Delta $ generated by the left translations of 
$L_{(p_i, {\tilde W}_j)}$ is  the semidirect product $\mathbb C^n \rtimes SU_{p_2}(n,\mathbb C)$  and the stabilizer 
 of $W_i$ in $\Delta $ is the semidirect product $\mathbb C^{p_i} \rtimes \Pi $, where $\Pi $ is an epimorphic 
image of the direct product $SU(p_1, \mathbb C) \times SU(p_2,\mathbb C) \times SO_2(\mathbb R)$.  

If  $F= \mathbb R$ then $\Delta $ is  the semidirect product  $\mathbb R^n  \rtimes SO_{p_2}(n,\mathbb R)^{\circ }$, where    
$SO_{p_2}(n,\mathbb R)^{\circ }$ is the connected component of $SO_{p_2}(n,\mathbb R)$, and the stabilizer 
of  $W_i$ in $\Delta $  is the semidirect product  $\mathbb R^{p_i} \rtimes (SO(p_1,\mathbb R) \times SO(p_2,\mathbb R))$. 
\end{Prop}
\Pro 
Clearly the topological manifold $L_{(p_i, {\tilde W}_j)}$ carries the differentiable 
structure of the real differentiable manifold $\Xi _{(p_i,{\tilde W}_j)}$ which is the 
topological  product of $T_{{\tilde W}_j}$ and  $\Sigma _{(p_1,p_2)}(F)$. 

According to Section $4$ the group $\Delta $ topologically generated by the left translations $\Xi _{(p_i,{\tilde W}_j)}$ is the semidirect product $\Delta=F^n \rtimes C$, where $C$ contains the commutator subgroup of $SU_{p_2}(n,F)$.

If $F=\mathbb C$ then $C=SU_{p_2}(n,\mathbb C)$ and the stabilizer $\hat{H}$ of $W_i$ in $\Delta $ is the semidirect product $T_{W_i} \rtimes \Phi$ with $\Phi =[U_{p_1}( \mathbb C) \times U_{p_2} (\mathbb C)] \cap  SU_{p_2}(n, \mathbb C)$ which is a maximal compact subgroup of $SU_{p_2}(n, \mathbb C)$ (\cite{Tits}, p. 28).  
The groups $SU_{p_2}(n, \mathbb C)$ and $\Phi $ are connected therefore the groups $\Delta $ and $\hat{H}$ are 
connected. Since  $\Delta $ is the topological product $\Xi_{(p_i, {\tilde W}_j)} \times {\hat H}= \Xi_{(p_i, {\tilde W}_j)} \times T_{W_i} \times \Phi$ it follows that  the manifold $\Xi _{(p_i,{\tilde W}_j)}$ and hence  the loop  
$L_{(p_i, {\tilde W}_j)}$ are   
diffeomorphic to an affine space. 

If $F=\mathbb R$ then $C$ is a subgroup of  $SO_{p_2}(n,\mathbb R)$ containing the commutator subgroup 
${\mathcal K}_{(n,p_2)}$. According to \cite{Dieudonne2} p. 57 the factor group 
$SO_{p_2}(n,\mathbb R)/ {\mathcal K}_{(n,p_2)}$ has order $2$. Hence  ${\mathcal K}_{(n,p_2)}$ is the connected 
component of $SO_{p_2}(n,\mathbb R)$. The group  $\Phi =[O_{p_1}(\mathbb R) \times O_{p_2}(\mathbb R)] \cap SO_{p_2}(n,\mathbb R)$ is not connected since the factor group $O(p_i, \mathbb R)/ SO(p_i, \mathbb R)$ has order $2$ (\cite{Porteous}, Corollary 9.37, p. 158) and  the product $\alpha _1 \alpha _2$ 
with $\alpha _i \in O(p_i, \mathbb R)$, but 
$\alpha _i \notin  SO(p_i, \mathbb R)$ for $i=1,2$, is an element of 
$SO_{p_2}(n,\mathbb R)$. The group  $SO_{p_2}(n,\mathbb R)$ is homeomorphic to the topological product 
$\Sigma _{(p_1,p_2)}(\mathbb R) \times \Phi$. Since  $SO_{p_2}(n,\mathbb R)$ has  two connected 
components and $\Phi $ is not connected the manifold $\Sigma _{(p_1,p_2)}(\mathbb R)$ is connected. It follows 
that the group $C$ generated by  $\Sigma _{(p_1,p_2)}(\mathbb R)$ is connected and hence isomorphic to the connected 
component $SO_{p_2}(n,\mathbb R)^{\circ }= {\mathcal K}_{(n,p_2)}$ of $SO_{p_2}(n,\mathbb R)$. Thus the group 
$\Delta =T \rtimes C$ is connected. Moreover $\Delta $ is the topological product  
$\Xi_{(p_i, {\tilde W}_j)} \times {\hat H}= \Xi_{(p_i, {\tilde W}_j)} \times T_{W_i} \times (\Phi \cap {\hat H})$. 
Since $\Delta $, $\Xi_{(p_i, {\tilde W}_j)}$ and $T_{W_i}$ are connected, the group $\Phi \cap {\hat H}$ is connected 
and hence a maximal compact subgroup of  $SO_{p_2}(n,\mathbb R)$. This yields that $\Xi_{(p_i, {\tilde W}_j)}$ and  
$L_{(p_i, {\tilde W}_j)}$ are diffeomorphic to an affine space. 

The group $\Delta $ is the topological product $\Xi_{(p_i, {\tilde W}_j)}  \times {\hat H}$. Thus  for the real dimension of $L_{(p_i, {\tilde W}_j)}$ one has 

\hspace{0.1cm}
\centerline{$\hbox{dim} L_{(p_i, {\tilde W}_j)}= \hbox{dim} \Xi_{(p_i, {\tilde W}_j)}=\hbox{dim} \Delta - \hbox{dim}  {\hat H}$}
\centerline{$\hbox{dim}_{\mathbb R} T_{W_i} +  \hbox{dim}_{\mathbb R} T_{{\tilde W}_j}+ \hbox{dim} SU_{p_2}(n,F) - \hbox{dim}_{\mathbb R} T_{W_i}- \hbox{dim}(C \cap {\hat H})$.}

\hspace{0.1cm}
If $F= \mathbb C$ then the group $\Phi =C \cap {\hat H}$ is an epimorphic image of the direct product $SU(p_1, \mathbb C) \times SU(p_2, \mathbb C) \times SO_2(\mathbb R)$ (cf. \cite{Tits}, p. 28). This yields 
\newline
\centerline{$ \hbox{dim} L_{(p_i, {\tilde W}_j)}=[(p_1+p_2)^2-1]+2 p_j-(p_1^2 -1) -(p_2^2-1) -1 = 2 p_j +2 p_1 p_2$}
 since the dimension of a unitary group  $SU_{k}(m,\mathbb C)$ is equal to $(m-1)^2+2(m-1)$ for $0 \le k \le m$ 
(\cite{Tits}, p. 26 and p. 28). 
It follows that $ L_{(p_i, {\tilde W}_j)}$ is diffeomorphic to  
$\mathbb R^{2(p_j+p_1 p_2)}$.  

The group $\Delta $ is the semidirect product $\Delta=\mathbb C^n \rtimes C$, where $C$ is the group  
$SU_{p_2}(n, \mathbb C)$ and the stabilizer ${\hat H}$ is the semidirect product $T_{W_i} \rtimes \Phi $, where 
 $\Phi $ is an epimorphic image of   
$SU(p_1, \mathbb C) \times SU(p_2, \mathbb C) \times SO_2(\mathbb R)$.

If $F=\mathbb R$ then $C \cap {\hat H}=SO(p_1,\mathbb R) \times SO(p_2,\mathbb R)$ (\cite{Tits}, p. 31 and p. 38). It follows that 
\newline
\centerline{$ \hbox{dim} L_{(p_i, {\tilde W}_j)}=\frac{1}{2} (p_1+p_2) (p_1+p_2-1)+p_j- \frac{1}{2} p_1 (p_1-1) -
 \frac{1}{2} p_2 (p_2-1)=$ }
\centerline{$p_j +p_1 p_2$.} 
Hence the loop  $L_{(p_i, {\tilde W}_j)}$ is diffeomorphic to 
$\mathbb R^{p_j +p_1 p_2}$.   

The group $\Delta $ is the semidirect product $\Delta=\mathbb R^n \rtimes C$, where $C$ is the group  
$SO_{p_2}(n,\mathbb R)^{\circ }$ and  the stabilizer ${\hat H}$ of $W_i$ in $\Delta $ is the semidirect product 
$\mathbb R^{p_i} \rtimes (SO(p_1, \mathbb R) \times SO(p_2,\mathbb R))$.

The loop $L_{(p_i, {\tilde W}_j)}$ is diffeomorphic to the manifold 
${\tilde W}_j \times {\cal Z}_i$ since ${\cal Z}_i$ is diffeomorphic to 
$\Sigma _{(p_1,p_2)}(\mathbb R)$. The mapping $(x, D^{\ast }) \mapsto x \vee D^{\ast }$ 
assigning to a point $x  \in {\cal A}_n=F^n$, $F \in \{\mathbb R, \mathbb C \}$ and to 
an element $D^{\ast }$ of the Gra\ss mannian manifold 
of the $(p_i-1)$-dimensional $F$-subspaces of the hyperplane $E$   the affine subspace $D$ containing $x$ and intersecting $E$ 
in $D^{\ast }$ is differentiable. Also the mapping $D \to D \cap {\tilde W}_j $ 
assigning to a $p_i$-dimensional affine $F$-subspace $D$ of ${\cal A}_n$ the point 
$D \cap {\tilde W}_j $ is differentiable. Since the loop realized on $\Sigma _{(p_1,p_2)}(F)$ 
by the multiplication $(2)$ is  differentiable, the  representation of $L_{(p_i, {\tilde W}_j)}$  on the manifold 
${\tilde W}_j \times {\cal Z}_i$ by 
the multiplication $(3)$ yields that  $L_{(p_i, {\tilde W}_j)}$ is differentiable.
\qed

\bigskip
\noindent
$\begin{array}{lcl}
\hbox{{\small \'Agota Figula}} & \quad & \hbox{{\small Karl Strambach}} \\
\hbox{{\small Mathematisches Institut }} & \quad &  \hbox{{\small Mathematisches Institut }} \\
\hbox{{\small der Universit\"at Erlangen-N\"urnberg}} & \quad & \hbox{{\small der Universit\"at 
Erlangen-N\"urnberg}} \\
\hbox{{\small  Bismarckstr. 1 $\frac{1}{2}$, }} & \quad & \hbox{{\small  Bismarckstr. 1 
$\frac{1}{2}$, }} \\ 
\hbox{{\small D-91054 Erlangen, Germany,}} & \quad & \hbox{{\small D-91054 Erlangen, Germany,}} \\
\hbox{{\small figula@mi.uni-erlangen.de}} & \quad & \hbox{{\small  strambach@mi.uni-erlangen.de}} \\
\hbox{{\small and}} & \quad &  \\
\hbox{{\small Institute of Mathematics,}} & \quad &  \\
\hbox{{\small University of Debrecen }} & \quad & \\
\hbox{{\small P.O.B. 12, H-4010 Debrecen,}} & \quad &  \\
\hbox{{\small Hungary, figula@math.klte.hu}} & \quad &
\end{array}$

\end{document}